\documentclass[11pt, a4paper]{article}
\usepackage[latin1]{inputenc}
\usepackage[T1]{fontenc}
\usepackage{amsmath}
\usepackage{amsfonts}
\usepackage{amssymb}
\usepackage{amsthm}
\usepackage[affil-it]{authblk}
\usepackage{enumitem}
\usepackage[normalem]{ulem}
\usepackage[mathscr]{euscript}
\usepackage{dsfont}
\usepackage{float}
\usepackage{times}
\usepackage{bbm}
\usepackage{xcolor}
\definecolor{LB}{rgb}{0.29, 0.63, 0.73}
\usepackage[colorlinks,linkcolor=blue,citecolor=blue,urlcolor=blue]{hyperref}
\usepackage[capitalise]{cleveref}
\usepackage[left=4cm, right=4cm, bottom=2cm, top=2cm]{geometry}

\usepackage[colorlinks,linkcolor=blue,citecolor=blue,urlcolor=blue]{hyperref}

\newtheorem{theorem}{Theorem}[section]
\newtheorem{theorem*}{Theorem}
\newtheorem{lemma}[theorem]{Lemma}

\theoremstyle{definition}
\newtheorem{definition}[theorem]{Definition}
\newtheorem{assumption}[theorem]{Assumption}

\theoremstyle{remark}
\newtheorem{remark}[theorem]{Remark}
\newtheorem{example}[theorem]{Example}

\numberwithin{equation}{section}

\crefname{example}{Example}{Examples}
\Crefname{example}{Example}{Examples}

\crefname{assumption}{Assumption}{Assumptions}
\Crefname{assumption}{Assumption}{Assumptions}

\crefname{condition}{Condition}{Conditions}
\Crefname{condition}{Condition}{Conditions}

\newcommand{\f}{\ensuremath{\frac}}

\newcommand{\C}{\ensuremath{\mathcal{C}}}
\newcommand{\F}{\ensuremath{\mathcal{F}}}

\renewcommand{\L}{\ensuremath{\mathcal{L}}}
\newcommand{\N}{\ensuremath{\mathbb{N}}}
\newcommand{\R}{\ensuremath{\mathbb{R}}}
\newcommand{\X}{\ensuremath{\mathcal{X}}}

\newcommand{\Y}{\ensuremath{\mathcal{Y}}}

\def\epsilon{\varepsilon}
\def\le{\leq}

\renewcommand{\P}{\ensuremath{\mathcal{P}}}

\newcommand{\define}{\ensuremath\triangleq}

\def\E{\mathbb E}
\def\EE{{\mathbf E}}

\renewcommand{\geq}{\geqslant}
\renewcommand{\leq}{\leqslant}
\def\CC{{\mathcal C}}

\def\${|\!|\!|}

\def\F{{\mathcal F}}
\newcommand{\vertiii}[1]{{\left\vert\kern-0.25ex\left\vert\kern-0.25ex\left\vert #1 \right\vert\kern-0.25ex\right\vert\kern-0.25ex\right\vert}}
\newcommand{\rom}[1]{(\textup{\uppercase\expandafter{\romannumeral#1}})}

\def\X{\mathbb{X}}
\def\CX{{\mathcal X}}

\def\X{{\mathbf X}}
\def\XX{{\mathbb X}}

\def\Y{{\mathbf Y}}
\def\YY{{\mathbb Y}}

\def\D{{\mathcal D}}
\def\CP{\mathcal P}
\def\err{\mathbf {Er}}
\def\FC{\mathscr{C}}
\def\CY{\mathcal Y}

\begin{document}
\title{ Rough Homogenisation with Fractional Dynamics  }
\author{Johann Gehringer and Xue-Mei Li \\Imperial College London \\ Dedicated to Sergio Albeverio on the occasion of his 80th birthday}
\maketitle

\begin{abstract}
We review recent developments of slow/fast stochastic differential equations, and also present a new result on
 Diffusion Homogenisation Theory with fractional and non-strong-mixing noise and providing new examples. 
 The emphasise of the review will be on the recently developed effective dynamic theory for two scale random systems with fractional noise: Stochastic Averaging and  `Rough Diffusion Homogenisation Theory'.  We also study the geometric models with perturbations to symmetries. 
\end{abstract}

\tableofcontents

\section{Introduction}


When we study the evolution of a variable / quantity, which we denote by $x$, we often encounter other interacting variables which either have the same scale as $x$ and are therefore treated equally, or are much smaller in size or slower in speed and are essentially negligible, or they might  evolve in a microscopic scale $\epsilon$, such variables are called the fast variables which we denote by $y$.  It  happens often that $y$ is approximately periodic or has  chaos behaviour or exhibits ergodic properties, then its effect on the $x$-variables can be analysed. During any finite time on the natural scale of $x$, the $y$-variable will have explored everywhere in its state space. As $\epsilon\to 0$,  the persistent effects from the fast variables will be encoded in the `averaged' slow motions through adiabatic transformation.   We then expect that  $\lim_{\epsilon\to 0} x_t^\epsilon$ exists; its limit will be autonomous, not depending on the $y$-variables.  
 In other words, the action of $y$  is transmitted adiabatically  to $x$ and the evolution of $x$ can be approximated by that of an autonomous system
called the effective dynamics of $x$. 

We explain this theory with  the two scale slow/fast  random evolution equation
\begin{equation}\label{eq:slow}
    dx_t^\varepsilon=f(x_t^\varepsilon,y_t^\varepsilon)\,dt+g(x_t^\varepsilon,y_t^\varepsilon)\,dB_t.
\end{equation}
Here $\epsilon>0$ is a small parameter,  $y_t^\epsilon$ is a fast oscillating noise, and $B_t$ is another noise. The stochastic processes $B_t$ and $y_t^\epsilon$  will be set on a standard probability space $(\Omega, \F, \P)$ with a filtration of $\sigma$-algebra $\F_t$.  Typically the sample paths of the stochastic processes $t\mapsto B_t(\omega)$, and $t\mapsto y^\epsilon_t(\omega)$ are not differentiable but have H\"older regularities. The equation is then interpreted as the integral equation:
\begin{equation}\label{eq:slow2}x_t^\varepsilon=x_0+ \int_0^t f(x_s^\varepsilon,y_s^\varepsilon)\,ds+\int_0^t g(x_s^\varepsilon,y_s^\varepsilon)\,dB_s.
\end{equation}
We take the initial values to be the same for all $\epsilon$. 
 The $x$-variables are usually referred to as the slow variables. 
If the fast dynamics depend on the slow variables, we refer this as  the `feedback dynamics'. If the fast variables do not depend on the slow variables, we have the `non-feedback dynamics'.   

If $(B_s)$ is a Brownian motion, the integral is an It\^o integral. The solutions of a It\^o stochastic differential equation (SDE) are Markov processes, they have continuous sample paths and therefore diffusion processes. Within the It\^o calculus realm,  the study of two scale systems began in the 1960s, almost as soon as a rigorous theory of It\^o stochastic differential equations was established, and has been under continuous exploration. 
In the averaging regime, the effective dynamics are obtained by averaging the original system in the $y$-variable. This  non-trivial dynamical theory is related to the Law of Large Numbers (LLN) and the ergodic theorems.  
The averaged dynamics for the Markovian system is expected to be again the Markov process whose Markov generator is obtained by averaging the $y$-components in the family of   Markov generators $\L^y$ of the slow variables with a parameter $y$.
Stochastic Averaging for Markovian ordinary differential equations was already studied in the 1960s and 1970s 
in \cite{Stratonovich61,Stratonovich63, Nelson,hasminskii68}, \cite{Borodin77, Freidlin-76, Skorohod,Kurtz70}. See also  \cite{Freidlin-76, Freidlin-Wentzell-03, Veretennikov,  Hairer-Pavliotis, Pavliotis-Stuart-Zygalakis, Skorohod-Hoppensteadt-Habib,Berglund-Gentz}.  Stochastic averaging with periodic and stationary vector fields  from the point of view of dynamical systems is a related classic topic, see \cite{Kryloff-Bogoliuboff}. For more recent work, see \cite{Neishtadt-91, Arnold-Imkeller-Wu, Bakhtin-Kifer,Kifer92}.  Stochastic averaging on manifolds for Markovian systems are studied in \cite{averaging, Ruffino, Li-OM-1, Catuogno-daSilva-Ruffino, Angst-Bailleul-Tardif,Birrell-Hottovy-Volpe-Wehr}.  
In the homogenisation regime, this theory is  linked to  Functional Central Limit Theorems (CLTs). In the classical setting this falls within the theory of diffusion creation, we therefore refer it as the diffusive homogenisation theory.
 A meta functional CLT  is as follows   \cite{Kipnis-Varadhan, Helland}: Let $f$ be such that $\int f d\pi=0$ where $\pi$ is  the unique invariant probability measure of a Markov  process $Y_s$, then $\f 1 {\sqrt t} \int_0^t f(Y_s) ds$ converges to a Markov process.  Such limit theorems are the foundation for diffusive homogenisation theorems and for studying weak interactions. If the dyanamics is Markovian, we naturally expect the limit of the slow variable to be another Markov process.  This theory is known as homogenisation; let us call this `diffusive homogenisation' to distinguish it from the settings where the dynamics is fractional. For Markovian systems, there are several well developed books, see e.g. \cite{Komorowski-Landim-Olla,Skorohod-Hoppensteadt-Habib}, see also \cite{Bensoussan-Lions-Papanicolaou11}.  Diffusive homogenisation was  studied in \cite{Kohler-Papanicolaou74,Hairer-Pardoux,Barret-vonRenesse,Dolgopyat-Kaloshin-Koralov}, see also  \cite{Liverani-Olla,Bakhtin-Kifer}.

An It\^o type stochastic differential equation is a good  model if the randomness is  obtained under the assumptions that there are a large number of independent (or  weakly correlated) components.
 However, long range dependence (LRD) is prevalent in mathematical modelling and observed in time series data such as economic cycles and data networks. One of the simplest LRD noise is the fractional noise, which is the `derivative' of fractional Brownian motions (fBM). Fractional Brownian motions was popularised by Mandelbrot and Van Ness \cite{Mandelbrot-VanNess} for modelling the long range dependent phenomenon observed by H. Hurst \cite{Hurst}. 
This is a natural process to use. Within the Gaussian class, those with stationary increments and the covariance structure
$\E(B_t-B_s)^2=(t-s)^{2H}$ are necessarily fractional Brownian motions with similarity exponent/Hurst parameter $H$. 
Fractional Gaussian fields and strongly correlated fields are used to  study critical phenomena in mathematical physics, see e.g. R. L. Dobrushin, G. Jona-Lasinio, G. Pavoliotti \cite{Dobrushin, Jona-Lasinio, Gallavotti-Jona-Lasinio},  and Ja. G. Sinai \cite{Sinai}.

 If $B_t$ is a fractional BM of Hurst parameter $H$,  the stochastic integral in (\ref{eq:slow2})  is a Riemann-Stieljes integral if $H>\f 	12$. Otherwise this can be understood in the sense of rough path integration or fractional calculus.  We explain the essence for this using the basics in the rough path theory. 
   The instrument for this is the  `Young bound'  \cite{Young36}:  if $F\in \CC^\alpha$ and $b\in \CC^\beta$ with $\alpha +\beta>1$  
  $$\left | \int_r^t (F_s-F_0) db_s \right| \lesssim  | F|_\alpha\, |b|_\beta \, (t-r)^{\alpha+\beta}.$$
With this, Young showed that the map $(F,b)\mapsto  \int _0^t F_s db_s$  is continuous where
$$\int_0^t F_s db_s= \lim_{|\CP|\to 0} \sum_{[u,v] \subset \CP}  F_u(b_u-b_v)\in \CC^\beta,$$
a Riemannan Stieljes integral / Young integral. 
This argument also established the continuity of the solutions to Young equations.   The solution $x_t$ inherits the regularity of  the driver
 and is in $C^{\beta}$.  This means if $H>\f 12$, the equation
   $$ dx_t=b(x_t)\,dt+F(x_t)\,dB_t$$
   below can be interpreted as a Young integral equation.
  It is well posed if  $F, b$ are in $BC^3$ \cite{Lyons94}. In general this type of SDE can be made sense of  if the Hurst parameter of the fBM $B_t$ is greater than $\f 14$, see \cite{Courtin-Qian, Decreusefond-Usutunel,Gradinaru-Nourdin-Russo-Vallois} and also \cite{Klingenhofer-Zahle, Mishura}.
 The study of stochastic evolution equations with fractional Brownian motions has since become popular, see eg. \cite{HuNual, Albeverio-Jorgensen-Paolucci,Nourdin-Nualart-Zintout-Rola,Flandoli-Gubinelli-Russo,Brezniak-Neerven-Salopek,  Cass-Hairer-Litteerer-Tindel,Garrido-Atienza-Schmalfuss} and\cite{Cont,Neuman-Rosenbaum} for their study in mathematical finance and in economics. See also \cite{Grahovac-Leonenko-Taqqu}.
 
 Despite of this popularity of stochastic equations with fractional noise,  there had not been much activity on the effect dynamics of a multi-scale systems.
A stochastic averaging with LRD fractional dynamics is obtained in
\cite{Hairer-Li} for the without feedback case and also for a feedback Markovian dynamics $y_t^\epsilon$ which solves the following interacting SDE:
\begin{equation}\label{fast}dy_t^\epsilon= \f 1 {\sqrt \epsilon}\sum_{k=1}^{m_2}  Y_k(x_t^\epsilon, y_t^\epsilon) \circ dW_t^k+
\f 1{\epsilon} Y_0(x_t^\epsilon, y_t^\epsilon)\,dt, \quad y_0^\epsilon =y_0
\end{equation}
where $W_t^i$ are independent Wiener processes. A uniform ellipticity is assumed of the equation.
 A priori the fast variables  $y_t^\epsilon \in \C^{\f 12 -}$. Since the sample paths of $B_t$ is in $C^{H-}$,  for the LRD case  where $H>\f 12$,
 (\ref{eq:slow}) is a Young equation.
More precisely,  for  $\omega$ fixed,  $$\Big( \int_0^t g(x_s^\epsilon, y_s^\epsilon)dB_s \Big) (\omega):=\int_0^t g(x_s^\epsilon(\omega), y_s^\epsilon(\omega))dB_s(\omega),$$   is a Riemannan Stieljes integral / Young integral.  The solution $x_t^\epsilon$ inherits the regularity of $B_t^H$ and is in $C^{H-}$.  These integrals are not defined with the classic It\^o stochastic calculus in any natural way,   it is therefore reasonable to use this pathwise interpretation.  A solution theory for the equation (\ref{eq:slow}) and (\ref{fast}) also exist, see  \cite{Guerra-Nualart, Hairer-Li}. 

As mentioned earlier, if $B_t$ is a BM,  the classic averaging theory states that the effective dynamics is the Markov process with its generator obtained by averaging the Markov generators of the slow variables.  This is obtained within the theory of  It\^o calculus. 
However  standard analysis within the integration theory does not lead to `pathwise'  estimates on $x_t^\epsilon$ that are uniform in $\epsilon$, which means a pathwise limit theory is not to be expected. It is clear that the `Young bounds' are totally ineffective for obtaining the essential uniform pathwise estimates for 
$\Big( \int_0^t g(x_s^\epsilon, y_s^\epsilon)dB_s \Big) (\omega)$.  When $\epsilon\to 0$ , the H\"older norm of the $y_\cdot^\epsilon$ is expected to blow up.  If $g$ does not depend on the fast variables,  it is of course possible to obtain pathwise bounds.
Indeed,  the sewing lemma of Gubinelli \cite{Gubinelli-lemma} and Feyel-de la Pradelle \cite{Feyel-delaPradelle}, which neatly encapsulates the main analytic estimates of both the work of Young  and that of Lyons \cite{Lyons94}, and has since become a fundamental tool in pathwise integration theory, does not provide the required estimates for the feedback dynamics.
Without any uniform pathwise estimates, the slow variables (for a generic equation) cannot be shown to converge for fixed fBM path.

%

In \cite{Hairer-Li},  a novel approximation, of the pathwise Young integrals by Wiener integrals, was introduced with the help of the stochastic sewing lemma of  L\^e \cite{Khoa}. This approach used, paradoxically,  the stochastic nature of the fractional Brownian motion in an essential way and, therefore, effectively departed the pathwise framework.
Since It\^o integrals and Wiener integrals are defined as an element of $L^2(\Omega)$,  the uniform estimates are  $L^p$-estimates and thus the limit theorem is an `annealed' limit.
It was shown,  \cite{Hairer-Li}, that $x_t^\epsilon$ converges in joint probability to the solution of the following equation with the same  initial data as $x_0^\epsilon$:
$$dx_t=\bar f(x_t)\, dt+\bar g(x_t) dB_t.$$
where $\bar f $ and $\bar g$ are obtained by directly averaging $f$ and $g$ respectively.
Stochastic averaging with fractional dynamics is  now a fast moving area, see \cite{Pei-Inaham-Xu, Pei-Inaham-Xu2, Rockner-Xie,bourguin2019typical} for more recent work see  \cite{Li-Sieber, Eichinger-Kuehn-Neamtu, bourguin2020discretetime}.

For the homogensation theory the main references are \cite{Gehringer-Li-fOU,Gehringer-Li-tagged}; see also \cite{Gehringer-Li-homo} which is the preliminary version of the previous two articles, equations of the form 
$$\dot x_t^\epsilon =h(x_t^\epsilon, y_t^\epsilon)$$
are studied. They can be used to  model the dynamics of a passive tracer in a fractional  and  turbulent random environment.

 In \cite{Gehringer-Li-fOU, Gehringer-Li-homo}, a functional limit theorem is obtained. The limit theorems are build upon the results in  \cite{Ustunel-Zakai, Nualart-Peccati,Nourdin-Nualart-Peccati, Nourdin-Nualart-Zintout, Bai-Taqqu, Pipiras-Taqqu,better-taqqu,Graphsnumber,Dobrushin-Major}.
In \cite{Gehringer-Li-tagged,Gehringer-Li-homo}, a homogenisation theorem for random ODE's for fractional Ornstein-Uhlenbeck processes.
See also \cite{Komorowski-Novikov-Ryzhik-14,Fannjiang-Komorowski-2000}.
Since the tools for diffusive homogenisation do not apply, we have to rely on a theorem from the theory of rough path differential equations. 
This approach is close to that in \cite{Kelly-Melbourne,   Catellier-Gubinelli, Chevyrev-Friz-Korepanov-Melbourne-Zhang,Perruchaud}, see also \cite{Friz-Gassiat-Lyons,Boufoussi-Ciprian,Al-Talibi-Hilbert}.  However, in these references only the dynamics are Markovian and the results are  of diffusive homogenisation type. In \cite{Gehringer-Li-tagged,Gehringer-Li-fOU, Gehringer-Li-homo}, the fast dynamics is a fractional Ornstein-Uhlenbeck process and the effective dynamics are not necessarily Markov processes and the limiting equation is a rough differential equation. 
We refer this theory as `rough creation' theory  and `rough homogenisation' theory.

 The study for  the stochastic averaging theory and the homogenisation theory for fractional dynamics has just started.
These theories departed from the classical theory both in terms of the methods of  averaging, the techniques, and the effective dynamics.  
We will  compare the methodologies and obtain the following intermediate result. 
for $y_t^\epsilon=y_{t/\epsilon}$ where $y_t$ is a stationary stochastic process  with  stationary distribution $\mu$.  Let $G_k$ be a collection of real valued  $L^2(\mu)$ functions on $\R$ we define
$$ X^{k,\epsilon}_t= \sqrt{\epsilon} \int_0^{\f t \epsilon} G_k(y_s) ds, \qquad \qquad X^{\epsilon}_t = \left( X^{1,\epsilon}_t, \dots, X^{n,\epsilon}_t \right)$$

\subsection{Main Results.} 
Our main result,   c.f. Section \ref{section1}, can be proved using the methods in \cite{Gehringer-Li-homo}. The statement is as follows:

\begin{theorem*}\label{theorem-1}
Let $y_t$ be a real valued stationary and ergodic process with stationary measure $\mu$.
Let  $G_k:\R\to \R$ be $L^2$ functions satisfying 
	\begin{equation}\label{conditioning}  \int_0^{\infty} \Vert \E \left[ G_k(y_s)  | \F_0 \right] \Vert_{L^2} \,ds < \infty.
	\end{equation} Suppose that furthermore the following moment bounds hold 
\begin{align*}
\Vert X^{k,\epsilon}_{s,t}  \Vert_{L^p} \lesssim \vert t-s \vert^{ \f  1 2}, \qquad \qquad 
\Vert \XX^{i,j,\epsilon}_{s,t} \Vert_{L^{\f p 2}} \lesssim \vert t-s \vert,
\end{align*}  
and the functional central limit theorem holds ( i.e. $X_t^\epsilon$ converges jointly in finite dimensional distributions to a Wiener process $X_t= \left(X^1_t, \dots , X^n_t \right)$).
Then the following statements hold.
\begin{enumerate}
\item  The canonical lift $\X^{\epsilon}=\left( X^{\epsilon}, \XX^{\epsilon} \right)$ converges weakly in $\FC^{\gamma}$ for $\gamma \in (\f 1 3, \f 1 2 - \f 1 p)$.
\item As $\epsilon\to 0$, the solutions of 
$$ \dot x_t ^\epsilon =\sum_{k=1}^N\f 1 {\sqrt \epsilon} \, f_k(x_t^\epsilon) \,G_k(y_t^\epsilon), \qquad x_0^\epsilon=x_0$$
converges  to the solution of the equation 
$ d{x}_t =\sum_{k=1}^n f_k(x_t) \circ d X^k_t$.
\end{enumerate}
\end{theorem*}
We will discuss examples of $y_t^\epsilon$ and $G_k$ for which the above holds including those discussed, see Section \ref{section2-4}.

\section{Homogenization via rough continuity}\label{section1}

With the Law of Large Number theory in place we now explain the homogenisation theory.
The homogenisation problem is about fluctuations from the average. We therefore take 
 \begin{equation}\label{2}   \dot x_t^\epsilon= h( x_t^\epsilon,y^{\epsilon}_t),
\end{equation}
where $h( x,y )$ is a function  averaging to zero and $h( x,y^{\epsilon}_t)$  is a fast oscillating fractional nature moving at a microscopic scale $\epsilon$. 
The aim is to obtain an `effective' closed equation whose solution $\bar x_t$ approximates $x_t^\epsilon$. This effective dynamics will have taken into accounts of the persistent averaging effects from the fast oscillations. Homogenisation for random ODEs has been dominated by the diffusion creation theory, with only a handful of exceptions where the limit is a fBM.  We will obtain a range of dynamics with local self-similar characteristics of the fractional Brownian motion, Brownian motion, and Hermite processes.

Recall that the  fractional Brownian motion is a Gaussian process with stationary increments. Its Hurst parameter $H$, given by its covariance structure: $\E(B_t-B_s)^2=(t-s)^{2H}$, indicates also the exponent in the power law decay of the corresponding fractional noise. Indeed, the correlation of two increments of the fractional Brownian motion of length $1$ and time $t$ apart, 
$$\E(B_{t+s+1}-B_{t+s})(B_{s+1}-B_s)=\f 12 (t+1)^{2H}+\f 12 (t-1)^{2H}-t^{2H}$$
which is approximately $2H t^{2H-2}$ for large $t$.

\subsection{CLT for stationary processes}

Given a stochastic process $y_s$ we are concerned with the question, whether for some scaling $\alpha(\epsilon)$ and function $G: \R \to \R$ the following term converges in the sense of finite dimensional distributions 
$$ X^{\epsilon}_t = \alpha(\epsilon) \int_0^{\f t \epsilon} G(y_s) ds. $$
Usual functional central limit theorems would set $\alpha(\epsilon) =  \sqrt{\epsilon}$ and the limit would be a Wiener process.
For Markovian noises these question has been studied a lot,  see eg \cite{Kipnis-Varadhan} and the book \cite{Komorowski-Landim-Olla} and the references therein. The basic idea is:   If $\L$ is the Markov generator
and the Poisson equation $\L=G$ is solvable, with solution in the domain of $\L$, then the central limit theorem holds. This follows from the martingale formulation for Markov processes.

Another kind of condition often imposed on the noise is some kind of mixing condition. Processes that satisfy these conditions are in a sense nicely behaved as they obey the usual CLT.
However, in this section we aim to treat cases in which non of this conditions are satisfied. For example the fractional Ornstein-Uhlenbeck process for $H> \f 1 2$ is neither Markovian nor obeys usual mixing assumptions. 
Looking at this question from a Gaussian perspective Rosenblatt \cite{Rosenblatt} gave an example of a stationary Gaussian sequence $X_k$ such that $ \frac  {1} {\sqrt{N} } \sum_{k=1}^N X_k $ does not converge, however the right scaling $\f  {1} {N^{\alpha} }
\sum_{k=1}^N X_k   $, for a suitable $\alpha > \f 1 2$ converges to the so called Rosenblatt process. Taqqu and Dobrushin \cite{Dobrushin} \cite{Taqqu-75,Taqqu} added to this work of so called non-central limit theorems. Philosophically, if the covariance function $\rho(j) = \E \left[ X_0 X_j \right]$ does not decay fast enough the limiting distribution can not have independent increments. The notions of short and long range dependence capture this idea. 
We say that a sequence is short range dependent if $ \sum_{j=1}^{\infty} \vert \rho(j) \vert < \infty$ and long range dependent otherwise.

We now want to discuss a method to conclude convergence to a Wiener process for these kind of processes.
In case of a Gaussian noise $y_s$ the rich toolboxes of Malliavin calculus, in particular the fourth moment theorem enables one to conclude limit theorems for a wide variety of situations, c.f. \cite{Nualart,Nourdin-Nualart-Zintout,Nourdin-Peccati,Nualart-Peccati}.
However, the method we are going to explore further relies  on a  martingale approximation method, see \cite{Jacod-Shiryaev} Thm 3.79.

The idea is to impose a condition on the functional instead on the noise in order to obtain a decomposition, as in the Markovian case, into a martingale and coboundary term.

Furthermore, the method of martingale approximation can be used to obtain convergence of the lifted process via  Theorem 2.2 in  \cite{Kurtz-Protter}, cf. \cite{roughflows,Kelly-Melbourne-16}.
Given a function $U$ satisfying Assumption \ref{assumption-conditioning} we may define the following $L^2$ martingale
$$
M_t = \int_0^{\infty} \E \left[ U(y_r) | \F_{t+} \right] - \E \left[ U(y_r) | \F_{0+} \right]dr.
$$
Now we may decompose $X^{\epsilon}$ as follows
$$ X^{\epsilon}_t =  \sqrt{\epsilon} \int_0^{ \f t \epsilon} U(y_r) dr = 
\sqrt{\epsilon} M_{\f t \epsilon} + \sqrt{\epsilon} \left( Z_{\f t \epsilon} - Z_0 \right),$$
where $ Z_t = \int_t^{\infty} \E \left[ U(y_r) | \F_{t+} \right] dr. $ As by assumption $\Vert Z_t \Vert_{L^2}$ is uniformly bounded we may drop the coboundary term and apply the Martingale central limit theorem.

\subsection{Fractional Ornstein Uhlenbeck as fast dynamics}

To illustrate the type of theorem we are seeking, we review the recently obtained result for the fractional Ornstein-Uhlenbeck process.
In   \cite{Gehringer-Li-tagged}, we 
studied  equation (\ref{multi-scale}) where $y_t^\epsilon$ is the fractional Ornstein-Uhlenbeck process. 
A fluctuation theorem from the average was obtained. We showed furthermore that the effective dynamics is the solution of (\ref{limit-eq}).

 \begin{definition}
A function $G\in L^2( \mu)$,   $G=\sum_{l=0}^\infty  c_l H_l$,
is said to satisfy the fast chaos decay condition with parameter $q\in \N$, if 
	$$\sum_{l=0}^\infty  {|c_l|}\; \sqrt{l!} \;(2q-1)^{\f l2}<\infty.$$ 
	 The lowest index $l$ with  $c_l\not =0$ is called the Hermite rank of $G$. If $m$ is the Hermite rank we define  $H^*(m)=m(H-1)+1$ .
 \end{definition}

 \begin{remark} \cite{Gehringer-Li-homo,Gehringer-Li-tagged}
 Let $y_t$ be the stationary fractional Ornstein-Uhlenbeck process with $H\in (0,1)\setminus\{ \f  1 2\}$. Then, for any real valued functions  $U \in L^2(\mu)$ with Hermite rank
 $m>0$  and $H^*(m)=m(H-1)+1<0$, we have
	\begin{equation}  \int_0^{\infty} \Vert \E \left[ U(y_s)  | \F_0 \right] \Vert_{L^2} \,ds < \infty.
	\end{equation}
Note that, without the conditioning on $\mathcal{F}_0$ the integral would be infinite due to the stationarity of $y_s$.
\end{remark}

Let $\alpha(\epsilon, H^*(m))$ be positive constants as follows, they depend on $m, H$ and $\epsilon$ and tend to $\infty$  as $\epsilon\to 0$,
 \begin{equation}\label{beta}
\begin{aligned}
&\alpha\left(\epsilon,H^*(m)\right) = \left\{\begin{array}{cl}
\f 1 {\sqrt{\epsilon}}, \, \quad  &\text{ if } \, H^*(m)< \f 1 2,\\
\f  1 {\sqrt{  \epsilon \vert \ln\left( \epsilon \right) \vert}}, \, \quad  &\text{ if } \, H^*(m)= \f 1 2, \\
\epsilon^{H^*(m)-1}, \quad  \, &\text{ if } \,  H^*(m) > \f 1 2.
\end{array}\right.
\end{aligned}
\end{equation}
 For $H^*(m) \not = \f 1 2$ this is equivalent to $\alpha(\epsilon,H^*(m)) =  \epsilon^{(H^*(m_k) \vee \f 1 2) -1} $.
 The following is proved in  \cite{Gehringer-Li-homo} for $H>\f 12$, see \cite{Gehringer-Li-tagged} for also the $H<\f 12$ case. 
\begin{theorem} \cite{Gehringer-Li-homo,Gehringer-Li-tagged}
Let $y_t$ be the fractional Ornstein-Uhlenbeck process with stationary measure $\mu$. Let $G_k$ be  $L^2\cap L^{p_k}$ functions with satisfies the fast chaos decay condition with parameter $q \geq 4 $ for $p_k$ sufficiently large (see blow).  We order the functions  $\{G_k\}$  so that their Hermite rank $m_k$ does not increase with $k$.
We also assume either $H^*(m_k) <0$ for $k \leq n$  or $H^*(m_k)> \f 1 2$ otherwise. Then, the solutions of 
	\begin{equation}\label{multi-scale}
\left\{ \begin{aligned} \dot x_t ^\epsilon &=\sum_{k=1}^N \epsilon^{(H^*(m_k) \vee \f 1 2) -1} \, f_k(x_t^\epsilon) \,G_k(y_t^\epsilon),\\
x_0^\epsilon&=x_0, \end{aligned}\right.  
\end{equation}
converges, as $\epsilon\to 0$  to the solution of the following equation with the same initial data:
\begin{equation}\label{limit-eq}
 d{x}_t =\sum_{k=1}^n f_k(x_t) \circ d X^k_t+\sum_{k=n+1}^N f_k(x_t) d X^k_t.
\end{equation}
Here $X^k_t$ is a Wiener process for $k \leq n$, and otherwise a Gaussian or a non-Gaussian Hermite process. The covariances between the processes are determined by the functions $G_k$,  for which there are explicit formulas. 
 In these equations, the symbol $\circ$ denotes the Stratonovich integral and the other integrals are in the sense of Young integrals.  
\end{theorem}

The conditions for $p_k$ are:
 \begin{assumption}
	\label{assumption-multi-scale}
	If $G_k$ has low Hermite rank,  assume $ H^*(m_k) - \f 1 {p_k} > \f 1 2$; otherwise assume $\f 1 2 - \f 1 p > \f 1 3$. Furthermore,
		\begin{equation}\label{Hoelder-sum>1}
		\min_{k\leq n} \left( \f 1 2 - \f 1 {p_k} \right) + \min_{n<k\leq N} \left( H^*(m_k) - \f 1  {p_k} \right)>1.
		\end{equation}
\end{assumption}

\subsection{The rough path topology}
To explain the methodology we first explain the necessities from the rough path theory.
If $X$ and $Y$ are H\"older continuous functions on $[0,T]$ with exponent $\alpha$ and $\beta$ respectively, such that $\alpha + \beta >1$, then by Young integration theory 
 $$\int_0^T Y dX=\lim_{\CP\to 0} \sum_{[u,v]\in \CP} Y_u(X_v-X_u).$$
  where $\CP$ denotes a partition of $[0,T]$. 
  Furthermore $(X,Y)\mapsto \int_0^T Y dX$ is a continuous map. Thus, for $X\in \C^{\f 12+}$ and $f\in \C_b^2$,
    one can make sense of a solution $Y$ to the Young integral equation $dY_s=f(Y_s) dX_s$. Furthermore the solution is  continuous with respect to both the driver $X$ and  the initial data, see \cite{Young36}.
If $X$ has H\"older continuity less or equal to $\f 12$,   this fails and one cannot define a pathwise integration for $\int X dX$ by the above Riemann sum anymore.  
Rough path theory provides us with a machinery to treat less regular functions by enhancing the process with a second order process, giving a better local approximation, which then can be used to enhance the Riemann sum and show it  converges.

 A rough path of regularity $ \alpha \in (\f 1 3 , \f 1 2)$,
is a pair of process $\X=(X_t, \XX_{s,t})$ where $(\XX_{s,t}) \in \R^{d \times d} $ is a two parameter stochastic processes
satisfying the following algebraic conditions: for $0\le s<u<t\le T$, 
$$\XX_{s,t}-\XX_{s,u}-\XX_{u,t}=X_{s,u} \otimes X_{u,t}, \qquad \qquad  \hbox{ (Chen's relation)} $$
where $X_{s,t}=X_t-X_s$, and $ (X_{s,u} \otimes X_{u,t})^{i,j}  = X^i_{s,u}  X^j_{u,t}$ as well as the following analytic conditions,
\begin{equation}\label{geo}
\Vert X_{s,t} \Vert \lesssim |t-s|^\alpha,  \qquad \Vert\XX_{s,t}\Vert \lesssim |t-s|^{2\alpha}.
\end{equation}
The set of such paths will be denoted by $\C^{\alpha}([0,T]; \R^d)$. The so called second order process $\XX_{s,t}$ can be viewed as a possible candidate for the iterated integral $\int_s^t X_{s,u} dX_u$. 

Given a  path $X$, which is regular enough to define its iterated integral, for example $X \in \C^1([0,T];\R^d)$, we define its natural rough path lift to be given by
$$\XX_{s,t}:=\int_s^t X_{s,u} dX_u.$$
It is now an easy exercise to verify that $\X = (X,\XX)$ satisfies the algebraic and analytic conditions (depending on the regularity of $X$), by which we mean Chen's relation and (\ref{geo}). Given two rough paths $\X$ and $ \Y$  we may define , for  $\alpha \in (\f 1 3, \f 1 2)$,  the following defines  a complete metric on $\C^{\alpha}([0,T]; \R^d)$, called the in-homogenous $\alpha$-H\"older rough path metric:
\begin{equation}\label{rough-distance}
\rho_\alpha(\X, \Y)=\sup_{s\not =t} \f {\Vert X_{s,t} -Y_{s,t} \Vert } {|t-s|^\alpha} 
+\sup_{s\not =t} \f {\Vert\XX_{s,t} -\YY_{s,t}\Vert }  {|t-s|^{2\alpha}} .
\end{equation}
We are also going to make use of the norm like object 
\begin{equation}
\Vert \X \Vert_{\alpha} = \sup_{s \not = t \in [0,T]} \f {\Vert X_{s,t}\Vert} {\vert t-s \vert^{\alpha} } +  \sup_{s \not = t \in [0,T]} \f {\Vert \mathbb{X}_{s,t}  \Vert^{\f 1 2}} {\vert t-s \vert^{\alpha}},
\end{equation}
where we denote for any two parameter process $\XX$ a semi-norm:
$$\|\XX\|_{2\alpha} := \sup_{s \not = t \in [0,T]} \f {\Vert \XX_{s,t}  \Vert} {\vert t-s \vert^{2\alpha}}.
$$

Our proof will be based on the following results:

\begin{lemma}\label{tightness-second-order}
	Let $\X^{\epsilon}$ be a sequence of rough paths with $\X(0)=0$ and $$\sup_{\epsilon \in (0,1]}\E \left(\Vert \X^{\epsilon}\Vert_{\gamma} \right)^{p} < \infty,$$
for some	$\gamma \in ( \f 1 3, \f 1 2 - \f 1 {p})$,
	then $\X^{\epsilon}$ is tight in $\C^{\gamma'}$ for every $\f 1 3 <\gamma' < \gamma$.
\end{lemma}

\begin{theorem}{\cite{Friz-Hairer}}\label{cty-rough}
	Let $Y_0 \in \R^m, \beta \in (\f1 3, 1), \, f \in \C^3_b(\R^m, {\mathbb L} (\R^d, \R^m)) $ and $\X \in \C^{\beta}([0,T],\R^d)$. Then, the differential equation
	\begin{equation}\label{example-sde}
	Y_t = Y_0 + \int_0^t f(Y_s) d\emph{X}_s 
	\end{equation}
	has a unique solution which belongs to $\mathcal{C}^{\beta}$. Furthermore, the solution map  $\Phi_f: ~\R^d\times \C^{\beta}([0,T], \R^d)
	\to  \D_{X}^{2\beta}([0,T],\R^m)$, where the first component is the initial condition and the second component the driver, is continuous.
	\end{theorem}

\subsection{Homogenization via rough continuity}
\label{section2-4}
\subsubsection{Main idea of the method}
Theorem \ref{cty-rough} has an interesting application to our homogenisation problem as weak convergence is preserved under continuous operations.
A simple equation for the demonstration is
$$
dx^{\epsilon}_t =\alpha(\epsilon)   f(x^{\epsilon}_t) G(y^{\epsilon}_t) dt,
$$
for a suitable choice of $\alpha$ and stochastic process $y^{\epsilon}_t$, we may rewrite this equation into a rough differential equation. To do so set 
$X^{\epsilon}_t = \alpha(\epsilon) \int_0^t G(y^{\epsilon}_s)\, ds$
thus,
$$
dx^{\epsilon}_t =  f(x^{\epsilon}_t) dX^{\epsilon}_t.
$$ 
To obtain a rough differential equation we now define the canonical lift of $X^{\epsilon}$ as  $\X^{\epsilon} = (X^{\epsilon}_t,\XX^{\epsilon}_{s,t})$, where $\XX^{\epsilon}_{s,t} = \int_{s}^t X^{\epsilon}_{s,r}\, dX^{\epsilon}_r$ and study
$$
dx^{\epsilon}_t =  f(x^{\epsilon}_t) d\X^{\epsilon}_t.
$$
As stated above to conclude weak convergence in a H\"older space of the solutions $x^{\epsilon}$ it is sufficient to obtain weak convergence of $\X^{\epsilon}$ in a rough path space $\FC^{\gamma}$  for some $\gamma \in ( \f 1 3, \f 1 2)$.
This can be done by following a two step approach.
Firstly, proving convergence in finite dimensional distributions for $\X^{\epsilon}$ and secondly tightness via moment bounds. See  \cite{Friz-Gassiat-Lyons,Kelly-Melbourne,roughflows,Gehringer-Li-tagged,Chevyrev-Friz-Korepanov-Melbourne-Zhang}.

Overall we have reduced the question to proving a functional central limit for $X^{\epsilon}_t = \alpha(\epsilon) \int_0^t G(y^{\epsilon}_s)\, ds$,  as well as for
$$\XX^{\epsilon}_{s,t} = \int_{s}^t X^{\epsilon}_{s,r} dX^{\epsilon}_r =\alpha(\epsilon)^2 \int_{s}^{t} \int_s^r G(y_r) G(y_u) du dr,$$
in a suitable path space. In one dimensions, by symmetry,  $\XX^{\epsilon}_{s,t} = \f {1}{2} \left( X^{\epsilon}_{s,t} \right)^2$, hence a continuous functional of $X^{\epsilon}$. This makes the one dimensional case quite simple.
However, when dealing with equations of the form
$$
dx^{\epsilon}_t =  \sum_{k=1}^N \alpha_k(\epsilon) f_k(x^{\epsilon}_t) G_k(y^{\epsilon}_t) dt,
$$
where the fast motions are channelled through functions of different scales, or 
$$
dx^{\epsilon}_t = \alpha(\epsilon) f(x^{\epsilon}_t,y^{\epsilon}_t),
$$
the canonical lift is more involved.

\subsubsection{Proof of Theorem \ref{theorem-1}} 
We apply the rough path theory and CLT theorem to diffusion creation.
To deal with components converging to  a Wiener process, to prove convergence of the second order process as mentioned above one usually relies on a martingale coboundary decomposition with suitable regularity, cf. \cite{Kelly-Melbourne,roughflows}.
Then, one can use the stability of weak convergence of the It\^o integral, Theorem 2.2 \cite{Kurtz-Protter}, for $L^2$ martingales. 

\begin{assumption}\label{assumption-conditioning}
	Given a stochastic process $y_t$ with stationary distribution $\mu$, a function $U$ in $L^2(\mu)$ is said to satisfy the `conditional memory loss' condition:
	\begin{equation} \int_0^{\infty} \Vert \E \left[ U(y_s)  | \F_0 \right] \Vert_{L^2} \,ds < \infty.
	\end{equation}
\end{assumption}
\begin{assumption}\label{moments}
A  vector valued process
$$ 
X^{\epsilon}_t = \left( X^{1,\epsilon}_t, \dots, X^{n,\epsilon}_t \right),
$$
is said to satisfy the rough moment condition if the following moment bounds hold 
\begin{align*}
\Vert X^{k,\epsilon}_{s,t}  \Vert_{L^p} &\lesssim \vert t-s \vert^{ \f  1 2}\\
\Vert \XX^{i,j,\epsilon}_{s,t} \Vert_{L^{\f p 2}} &\lesssim \vert t-s \vert.
\end{align*}  
\end{assumption}

Given a stationary and ergodic process $y_t$ and functions $G_k$, we define 
$ X^{k,\epsilon}_t= \sqrt{\epsilon} \int_0^{\f t \epsilon} G_k(y_s) ds$.  Suppose that there exists a Wiener process $X_t= \left(X^1_t, \dots , X^n_t \right)$ such that 
$ 
X^{\epsilon}_t \to X_t
$
 in finite dimensional distributions. 
 To prove Theorem \ref{theorem-1}, by the continuity theorem, it is sufficient to show that the canonical lift $\X^{\epsilon}=\left( X^{\epsilon}, \XX^{\epsilon} \right)$ converges weakly in $\FC^{\gamma}$ for $\gamma \in (\f 1 3, \f 1 2 - \f 1 p)$.

 \begin{remark}
In case the assumptions of  Theorem \ref{theorem-1} are satisfied one has in particular, 
$$\E [ X^j_t X^l_s ] = 2 ( t \wedge s ) \int_0^{\infty} \E [ G_j(y_r) G_l(y_0) ] dr.$$
\end{remark}

{\bf Proof of Theorem \ref{theorem-1}.}
To prove Theorem \ref{theorem-1} one may argue similarly as in section 3.3 in \cite{Gehringer-Li-tagged}. 
Recall that we assume that $y_t$ and $G_k$  satisfy the conditional memory  loss condition in Assumption \ref{assumption-conditioning}, 
and $X^\epsilon$ satisfies Assumption \ref{moments}.

Firstly, due to Assumption \ref{assumption-conditioning} we may decompose each $X^{k,\epsilon}$ as follows,
$$ X^{k,\epsilon}_t =  \sqrt{\epsilon} \int_0^{ \f t \epsilon} G_k(y_r) dr = 
\sqrt{\epsilon} M^k_{\f t \epsilon} + \sqrt{\epsilon} \left( Z^k_{\f t \epsilon} - Z^k_0 \right),$$
where $ Z^k_t = \int_t^{\infty} \E \left[ G_k(y_r) | \F_{t+} \right] dr $ and $
M^k_t = \int_0^{\infty} \E \left[ G_k(y_r) | \F_{t+} \right] - \E \left[ G_k(y_r) | \F_{0+} \right]dr.
$
By the construction, $M^k_t$ is a martingale and $ \Vert Z^k_t \Vert_{L^2}$ is bounded by Assumption \ref{assumption-conditioning}. Hence, the term $\sqrt{\epsilon} \left( Z^k_{\f t \epsilon} - Z^k_0 \right) $ converges to $0$ in $L^2$.
Thus, also the multidimensional martingales $$(\sqrt{\epsilon} M^1_{\f t \epsilon} , \dots, \sqrt{\epsilon} M^n_{\f t \epsilon})$$ 
converge jointly to the $n$-dim Wiener process,
$$
X_t=(  X^1_t, \dots,  X^n_t).
$$
When computing the iterated integrals $\XX^{i,j,\epsilon}_{s,t} = \int_s^t X^{i,\epsilon}_{s,r} dX^{j,\epsilon}_r$, we obtain for the diagonal entries $\XX^{i,i,\epsilon}_{s,t} = \f  1 2  \left( X^{i,\epsilon}_{s,t} \right)^2$ by symmetry as above. However, in case $i \not = j $ we need to argue differently. Using the same martingale coboundary decomposition as above one can show the following lemma.
\begin{lemma}[See the Appendix]\label{lemma-in-appendix}
For $L(\epsilon)= \lfloor \f t \epsilon \rfloor$,
$$ \epsilon \int_0^{ \f t \epsilon} \int_0^s G_i(y_s) G_j(y_r)dr = \epsilon \sum_{k=0}^{ L(\epsilon)} ( M^i_{k+1} - M^i_k ) M^j_k + tA^{i,j} + \err(\epsilon),$$
where $A^{i,j} = \int_0^{\infty} G_i(y_s) G_j(y_0) ds$ and $\err(\epsilon)$ converges to $0$ in probability.
\end{lemma}
Now, defining the cadlag martingales $M^{i,\epsilon}_t = \sqrt{\epsilon} M^i_{ \lfloor \f t \epsilon \rfloor }$ one may identify the sum above as It\^o integral. Furthermore, by Assumption \ref{assumption-conditioning} these martingales are bounded in $L^2$, hence, using Theorem 2.2 in \cite{Kurtz-Protter} we obtain,
$$ \left(M^{k,\epsilon}_t, \int_0^t M^{i,\epsilon}_s dM^{j,\epsilon}_s \right) \to \left( X_t, \XX_{0,t} \right),$$
where $\XX$ is given by It\^o integrals, in the sense of finite dimensional distributions. Additionally, the moment bounds guarantee that the convergence actually takes place $\FC^{\gamma}$ for $\gamma \in (\f 1 3 , \f 1 2 - \f 1 p)$.\\
\\
An example that satisfies above conditions is given by choosing $y_t$ as the fractional Ornstein Uhlenbeck process and functions $G_k$ that satisfy,

\begin{assumption}
	Each $G_k $ belongs to $L^{p_k}(\mu)$, where $p_k >2$, and has Hermite rank $m_k\ge 1$. Furthermore,
	\begin{enumerate}
		\item [(1)] Each $G_k$ satisfies the fast chaos decay condition with parameter $q \geq 4 $.
		\item[(2)] For each $k$,  $\f 1 2 - \f 1 {p_k} > \f 1 3$ and $H^*(m_k) <0.$
	\end{enumerate}
\end{assumption}
where $H^*(m) = (H-1)m +1$, the Hermite rank is as defined in section \ref{subsection-non-wiener-homo} .

\subsection{Examples satisfying Assumption \ref{assumption-conditioning}}

\subsubsection{Strong Mixing environment}
Using this method one can also treat the classical setup. Given a stationary stochastic process $y_t$ and assume it is strong mixing with mixing rate $\alpha(s)$.
$$
\alpha(t) = \sup_{A \in \F^0_{\infty}, B \in \F_{t}^{\infty} } \vert \P(A \cap B) - \P(A) \P(B) \vert. 
$$
By Lemma 3.102 in \cite{Jacod-Shiryaev}, given $G : \R \to \R$ such that $\Vert G(y_0) \Vert_{L^1} < \infty$ the following inequality holds,
$$ \Vert \E\left[ G(y_t) | \F_0\right] - \E\left[ G(y_t) \right] \Vert_{L^q} \le  2\left(2^{\f 1 q} +1\right)  \alpha(t)^{\f 1 q - \f 1 r} \Vert G(y_t) \Vert_{L^r} , 
$$
where $ 1 \leq q \leq r \leq \infty$.
\begin{lemma}
Let  $y_t$ be a stationary process with mixing rate $\alpha(t)$ and $G \in L^r(\mu)$  is centred.
Then Assumption \ref{assumption-conditioning} holds if 
$$\int_0^\infty \alpha(t) ^{\f 12-\f 1r} dt<\infty.$$
\end{lemma}

\subsubsection{Volterra kernel moving averages}
Let $B$ denote a fBm of Hurst parameter $H$ and set 
$y_t = \int_{-\infty}^t K(t-s) dB_s$, where $K$ denotes a  kernel $K$ such that $\Vert y_t \Vert_{L^2}=1$.
Using the decomposition $B_t-B_k=\tilde B_t^k+\bar B_t^k$, where
$ \bar B_t^k  = \int_{-\infty}^k (t-k)^{H - \f 1 2 } - (k-r)^{H- \f 1 2} dr$ and $\tilde B_t^k = \int_k^t (t-r)^{H- \f 1 2} dW_r$, we may also decompose 
\begin{align*}
y_t &= \int_{-\infty}^{t} K(t-r) dB_r
= \int_{-\infty}^k K(t-r) dB_r +   \int_{k}^{t} K(t-r) d( B_r - B_k)\\
&=  \left( \int_{-\infty}^k K(t-r) dB_r +\int_{k}^{t} K(t-r) d\overline{B}^k_r \right)+  \int_{k}^{t} K(t-r)
d\tilde{B}^k_r\\
&= \overline{y}^k_t + \tilde{y}^k_t.
\end{align*}
It was shown in \cite{Hairer05} that the term $\tilde{B}^k_t$ is independent of  $\mathcal{F}_k$ and $\overline{B}^k_t$ is  $\mathcal{F}_k$ measurable, hence  $\tilde{y}^k_t$ is independent of  $\mathcal{F}_k$  and $\overline{y}^k_t$ is  $\mathcal{F}_k$ measurable. Moreover, both terms are Gaussians. We set $y^k_t = y_t$ for $k \geq t$.
Using an expansion into Hermite polynomials one can show the following:

\begin{lemma} Let $y_t = \int_{-\infty}^t K(t-s) dB_s$ such that $\Vert y_t \Vert_{L^2} =1$ be given. If the kernel $K$ is such that
\begin{align*}
\int_{k-1}^{\infty}  \left( \E  \left[ {\left(\overline{y}^k_s\right)}^2 \right] \right)^{\f m 2}  dt,
\end{align*}
then Assumption \ref{assumption-conditioning} hold.
\end{lemma}
An example of this is the fractional Ornstein-Uhlenbeck process.

\subsubsection{Effective dynamics driven by  non-Wiener process}\label{subsection-non-wiener-homo}
As mentioned above for Gaussian noises one may obtain a variety of limits. The Hermite rank plays a central role in the analysis of limit theorems with Gaussian noises. Let $\gamma$ denote the standard normal distribution and $y_t$ a stationary process with distribution $\gamma$, then each $G \in L^2(\gamma)$ admits an expansion into Hermite polynomials, $G(y_s)=\sum_{k=0}^{\infty} c_k H_k(y_s)$. The Hermite rank is now defined as the rank of the smallest non zero Hermite polynomial, thus, a function with Hermite rank $m$ can be written as $G(y_s) =  \sum_{k=m}^{\infty} c_k H_k(y_s)$ with $c_m \not = 0$.  Hermite polynomials are mutually orthogonal and satisfy 
$\E\left( H_k(y_s) H_j(y_r) \right) = \delta(k-j)\rho( \vert s-r \vert)^k $, where $\rho$ denotes the correlation function of $y_s$. In case $\rho$ admits only an algebraic decay, ie. the fractional OU process satisfies $\rho(t) \lesssim 1 \wedge t^{2H-2}$, higher order polynomials accelerate the correlation decay. Thus, the picture one obtains is that there is a critical $m$ such that $H_k(y_s)$ is short range dependent if $ k \geq m$ and long range dependent otherwise. In case $ (H-1)k+1 > \f 1 2 $ the process is long range dependent and short range dependent if $(H-1)k+1 < \f 1 2 $, the border line case is long range dependent as well, however the sum of the correlations only diverges logarithmically, which leads to a separate behaviour which will not be discussed here.
Thus, if a function $G$ has Hermite rank $m$ such that $(H-1)k+1 > \f 1 2$  one obtains convergence to Hermite processes for $\alpha(\epsilon) = \epsilon^{(H-1)m+1}$. As this processes have H\"older regularity greater than $\f 1 2$ they can be treated within the Young framework.

\begin{remark}
 To combine both frameworks in \cite{Gehringer-Li-tagged} the assumption was made that the H\"older regularity of Wiener components plus the ones for Hermite components is bigger than $1$, hence the joint lifts are irrelevant in the limit and well defined as Young integrals.
\end{remark}

 \section{Recent progress  on  slow/fast Markovian Dynamics }

In order to compare the methods, we will explain the classical Stochastic Averaging Theory and the Diffusion Homogenisation theory with Markovian Dynamics, which have been continuously re-inventing itself since  the 1960s and 1970s. They are typically a system of two-scale stochastic equations as follows
\begin{equation}\label{sdes}
\left\{\begin{aligned} dx_t^\epsilon&=\sum_{k=1}^{m_1} X_k(x_t^\epsilon, y_t^\epsilon)\circ d\tilde W_t^k+
X_0(x_t^\epsilon, y_t^\epsilon) \,dt, \quad x_0^\epsilon=x_0;\\
dy_t^\epsilon&= \f 1 {\sqrt \epsilon}\sum_{k=1}^{m_2}  Y_k(x_t^\epsilon, y_t^\epsilon) \circ dW_t^k+
\f 1{\epsilon} Y_0(x_t^\epsilon, y_t^\epsilon)\,dt, \quad y_0^\epsilon =y_0.\\
\end{aligned}\right.
\end{equation}
where  $W_t^i, \tilde W_t^i$ are independent Brownian motions. 
Such models have the flavour of the following multi-scale system in dynamical system
$$\dot x_t^\epsilon=f(x_t^\epsilon, y_t^\epsilon), \qquad \qquad \dot y_t^\epsilon =\f 1 \epsilon Y_0(x_t^\epsilon, y_t^\epsilon).$$
In dynamical system, the fast dynamics is often periodic or has chaotic behaviour.

Because these are well known,  we will focus on the newer developments that applies to non-linear state spaces which is sufficient to illustrate the underlying 
ideas and the differences between the theories within the Markovian dynamics and that with the fractional dynamics. 
Since the scale in a multi-scale system  are note  naturally separated, we will need to use geometric methods to separate them, the slow and fast variables so obtained often lives in a non-linear space. For example  if we take an approximately integrable Hamiltonian system the natural state space is a torus. If the 
slow-fast system is obtained by using symmetries, the state space of the fast motions  are Lie groups. SDEs with symmetries are popular topics, see \cite{DFMU,albeverio-DeV-Morando-Ugolini, Elworthy-LeJan-Li-book-2, Takao, Elworthy-Li-Intertwining, ELL-equivariant}. See also  \cite{conservation,  Li-OM-1} for perturbation to symmetries.

 \subsection{The Basic Averaging Principal}
For $y\in \R^n$, let $\sigma_i(\cdot, y) :\R^d\to \R^d$, and for $x\in \R^d$,  $ Y_i(x, \cdot): \R^n\to \R^n$. 
The averaging principle states that if $(x_t^\epsilon, y_t^\epsilon)$ is the solution of the following stochastic differential equation, with It\^o integrals, 
\begin{equation}\label{sdes-2}
\left\{\begin{aligned} dx_t^\epsilon&=\sum_{k=1}^{m_1}  \sigma_k(x_t^\epsilon, y_t^\epsilon) d\tilde W_t^k+
\sigma_0(x_t^\epsilon, y_t^\epsilon) \,dt, \quad x_0^\epsilon=x_0;\\
dy_t^\epsilon&= \f 1 {\sqrt \epsilon}\sum_{k=1}^{m_2}  Y_k(x_t^\epsilon, y_t^\epsilon) dW_t^k+
\f 1{\epsilon} Y_0(x_t^\epsilon, y_t^\epsilon)\,dt, \quad y_0^\epsilon =y_0,\\
\end{aligned}\right.
\end{equation}
then $x_t^\epsilon$ converges to a Markov process.  
This principle requires technical verification.  

The popular format and assumptions for the Stochastic Averaging Principle is as follows, see \cite{hasminskii68} and \cite{Veretennikov}. 
Let $Y_s^x$ be the solution to the equation below with frozen slow variable (we assume sufficient regularity assumption so that there equation has a unique strong solution,)
$$dy_t^x=\sum_{k=1}^{m_2}  Y_k(x, y_t^x)  dW_t^k+
 Y_0(x_, y_t^x)\,dt, \quad y_0^x =y_0.$$
 suppose that there exists a unique invariant probability measure $\mu^x$.
   Assume furthermore that the coefficients of the equations are Lipschitz continuous in both variables and suppose that 
  there exist functions  $\bar a_{i,j}$ and $\bar b$ such that  on $[0,T]$,
\begin{equation}\label{ergodic-assumption}
\begin{aligned} 
&\left|\f 1 t \EE \int_0^t b(x, Y_s^x)ds-\bar b(x)\right|\lesssim(|x|^2+|y_0|^2+1),\\
&\left|\f 1 t \EE \int_0^t \sum_{k}\sigma_k^i\sigma_k^j(x, Y_s^x)ds-\bar a_{i,j}(x)\right|\lesssim (|x|^2+|y_0|^2+1).
\end{aligned}
\end{equation}  
Then, $x_t^\epsilon$ converges weakly to the Markov process with generator 
$$\bar \L =\f 12 \bar a_{i,j} (x) \f {\partial^2}{\partial x_i\partial x_j}+\bar b_k(x)\f {\partial}{\partial x_k}.$$
The notation $\bar h$  denotes the average of a function $h$, i.e.  $\bar h(x)=\int h(y) \mu^x(dy)$.
In \cite{hasminskii68},  boundedness in $y$ is assumed in (\ref{ergodic-assumption}).
This is replaced by the quadratic growth in \cite{Veretennikov}, where a uniform ellipticity is also used to replace the regularity assumption needed on $x\mapsto \mu^x$.
   
  In the next section we  outlined conditions posed directly on the coefficients of the equations, under which these assumptions hold.

\subsection{Quantitative locally uniform  LLN }
The law of large numbers is the foundation for stochastic averaging theory.
Let $y_t^\epsilon$ be an ergodic Markov process on a state space $\CY$ with invariant measure $\pi$, then the Birkhoff's ergodic theorem holds, i.e.
$$\lim_{t\to \infty}\left|\f 1t \int_0^t h(x_s)ds-\int h(y)\pi(dy) \right|\to 0.$$
  If $\L$ is a strictly elliptic operator, this is seen easily by the martingale inequalities and Schauder's  estimates.  There is typically  a rate of convergence of the order $\f 1 {\sqrt t}$. Such results are classic for elliptic diffusions on $\R^n$. 
For Brownian motion on a compact manifolds this is proven in \cite{Ledrappier}, see also \cite{Enriquez-Franchi-LeJan}, we generalise this result to non-elliptic operators and obtain a quantitative estimate.  
Furthermore, we obtain locally uniform estimates for diffusion operators $\L_x$ where  $x\in \CX$ is a parameter.  Our main contribution  is to obtain quantitative estimates that are  locally uniform in $x$.
We indicate one such result below.
\begin{definition}\label{def-hormander}
Let $X_0, X_1, \dots, X_k$ be smooth vector fields.
\begin{enumerate}
\item The differential  operator $\sum_{k=1}^m (X_i)^2+X_0$  is said to satisfy {\it H\"ormander's condition}  if  $\{X_k, k=0, 1, \dots, m\}$ and their iterated Lie brackets generate the tangent space at each point.
\item The differential  operator $\sum_{k=1}^m (X_i)^2+X_0$, is said to satisfy {\it strong H\"ormander's condition}  if  $\{X_k, k=1, \dots, m\}$ and their iterated Lie brackets generate the tangent space at each point.

\end{enumerate}
\end{definition}

Let $s\ge 0$, let $dx$ denote the volume measure of a Riemannian manifold $G$ and let $\Delta$ denote the Laplacian. 
 If $f$ is a $C^\infty$ function we define its Sobolev norm to be
 $$\|f\|_{s}= \left( \int_G f(x)(I+\Delta)^{s/2}f(x) \,dx\right)^{\f 12}$$

If the Strong H\"ormander's condition holds and $\CY$ is compact, then the Markov process with generator  $\sum_{k=1}^m (X_i)^2+X_0$ has a unique invariant probability measure.
We state a theorem for $\CY$ compact, a version with $\CY$ not compact can also be found in \cite{conservation}.
\begin{theorem}\label{quantitative}
\cite{conservation}
\label{lln}
Let $\CY$ be a compact manifold. Suppose that $Y_i(x,\cdot)$ are bounded  vector fields on $\CY$ with bounded derivatives and $C^\infty$ in both variables. Suppose that for each
$x\in \CX$, $$\L_x=\f 12\sum_{i=1}^m  Y^2_i(x, \cdot) +Y_0(x, \cdot)$$ satisfies H\"ormander's condition
and has a unique invariant probability measure which we denote by $\mu_x$. Then the following statements hold .
\begin{enumerate}
\item[(a)]  
  $x\mapsto \mu_x$ is locally Lipschitz continuous in the total variation norm.
\item[(b)] 
  For every  $s>1+\f{\dim(\CY)}{2}$ there exists a positive constant $C(x)$, depending continuously in $x$,  such that for every smooth  function $f:\CY \to \R$, 
\begin{equation}
\label{llln-1}
\left| {1\over T} \int_t^{t+T} f(z_r^x) \;dr- \int_G f(y) \mu_x(dy) \right|_{L_2(\Omega)} 
\le C(x)\|f\|_s\f 1{\sqrt T},
\end{equation}
where  $z_r$ denotes an $\L_x$-diffusion.
\end{enumerate}
\end{theorem}
The proof for this follows from an application of It\^o's formula, applied to the solution of the Poisson equation $\L_x h=f(x,\cdot)$ where  $\int f(x, y) \mu_x(dy)=0$. For such functions,
$${1\over T} \int_0^T f(x, z_r^x)dr
={1\over T} \left(g(x, z_T^x)-g(x, y_0)\right)
-{1\over T}\left(\sum_{k=1}^{m_2}  \int_0^T dg(x, \cdot)(Y_k(x, z_r^x) ) dW_r^k\right).$$
(We take $t=0$ for simplicity.) It then remain to bound  the  supremum norm of
$dg(x, \cdot)$, which is a consequence of the sub-elliptic estimates of H\"ormander.
%

  \subsection{Averaging with H\"ormander's conditions} 
  Let $\CX, \CY$ be smooth manifolds.
We take a family of  vector fields $\sigma_i(\cdot, y)$ on $\CX$  indexed by $y\in \CY$ and a family of vector fields $Y_i(x,\cdot)$ with parameter $x\in \CX$.
The vector field $\sigma_i(\cdot, y)$ acts on a  real function  $h:\CX\to \R$ so that $\sigma_i(\cdot, y)h$ is the derivative of $h$ in the direction of $\sigma_i(\cdot, y)$.
If $X$ is a vector field we denote by $Dh(X(x))$ the derivative of the function $h$ in the direction of $X(x)$ at the point $x$. The function so obtained is also denoted by $\L_Xh $, or $Xh$, or $Dh(X)$.

 The assumption (\ref{ergodic-assumption}) can in fact be verified with ergodicity and regularity conditions on the coefficients.
We state such a result on manifolds. We switch to Stranovich integrals and denote this by $\circ$. We also denote by $\tilde \sigma_0$ and $\tilde Y_0$ the effective drifts (including the Stratonovich corrections).  Consider  
\begin{equation}\label{sdes-2}
\left\{\begin{aligned} dx_t^\epsilon&=\sum_{k=1}^{m_1}  \sigma_k(x_t^\epsilon, y_t^\epsilon) \circ d\tilde W_t^k+
\sigma_0(x_t^\epsilon, y_t^\epsilon) \,dt, \quad x_0^\epsilon=x_0;\\
dy_t^\epsilon&= \f 1 {\sqrt \epsilon}\sum_{k=1}^{m_2}  Y_k(x_t^\epsilon, y_t^\epsilon) \circ dW_t^k+
\f 1{\epsilon} Y_0(x_t^\epsilon, y_t^\epsilon)\,dt, \quad y_0^\epsilon =y_0.
\end{aligned}\right.
\end{equation}

 Suppose that, for each $x$,    $\L_x=\f 12\sum_{i=1}^m  Y^2_i(x, \cdot) +Y_0(x, \cdot)$ has a unique invariant probability measure which we denote by $\mu^x$. 
 For $h:\CX\to \R$  a smooth function with compact support, we define
$$\bar \L h(x)=\f 12   \sum_k  \int_\CY  \Big(D^2h)_x (\sigma_k(x, y), \sigma_k(x, y)) \Big) (y)  \mu^x(dy) +\int _\CY (Dh)_x\Big( \sigma_0 (x, y )\Big)  \mu^x(dy).$$
The key condition for the convergence of the slow motion  is the H\"ormander's condition. More precisely we will make use of the conclusion of Theorem \ref{quantitative}, from which
 we obtain the following estimates  on a compact subset $D$ of $\CX$:
 \begin{equation} \E\sum_{i=0}^{N-1}\left|  \int_{t_i}^{t_{i+1}}f\left(x_{t_i}^\epsilon, y^{x_{t_i }^\epsilon}_{r}\right) \,ds- \Delta t_i\, f\left( x_{t_i}^\epsilon\right)  \right|
\le c\, T  \,\lambda(\f {\Delta t_i} \epsilon)\sup_{x\in D}  \,\left \|f(x, \cdot)-\bar f(x)\right\|_{s}.
\end{equation}
Here $\lambda(t)$ is a function converging to zero as $t\to \infty$. With this we obtain the dynamical Law of Large Number as follows.

\begin{theorem} \cite{conservation}\label{LLN}
Let $\CY$ be compact and $Y_i$ are in $C^2\cap BC^1$. Suppose that $\tilde \sigma_0$  and $\sigma_i$ are $C^1$, where $i=1,\dots, m$,
and  $\L_x$ satisfies  H\"ormander's condition and  that it  has a unique invariant probability measure.
Suppose that {\bf one} of the following  two statements holds.
\begin{enumerate}
\item [(i)] Let $\rho$ denote the distance function on $\CY$. Suppose that  $\rho^2$ is smooth and  $$\f 12 \sum_{i=1}^m \nabla d \rho^2(\sigma_i(\cdot,y), \sigma_i(\cdot, y))+ d\rho^2 (\tilde \sigma_0(\cdot,y ))\le  K(1+\rho^2(\cdot)), \quad \forall y\in G.$$

 \item [(ii)]  The sectional curvature of $\CY$ is bounded. There exists a constant  $K$ such that  $$\sum_{i=1}^m| \sigma_i(x,y) |^2\le K(1+ \rho(x)), \quad\left  |\sigma_0(x,y) \right| \le  K(1+ \rho(x)), \quad \forall x\in N, \forall y\in G.$$

\end{enumerate}
Then  $\{x_t^\epsilon, \epsilon>0\}$ converges weakly,  on any compact time intervals, to the Markov process with the Markov generator $\bar \L$.
\end{theorem}

This theorem is useful when we have an stochastic differential equation on a manifold which is invariant under a group action, and when we study its small perturbations.

\subsubsection{Geometric Models}
The  LLN  can be used to study the following models.
\begin{example}
\begin{enumerate}
\item [(1)] {\bf Approximately  Integrable Hamiltonian systems. }
If  $H: M\to \R$  is a smooth Hamiltoaian function on a symplectic manifold or on $\R^{2n}$ with its standard symplectic structure, we use $X_{H}$ to denote its Hamiltonian vector fields. Let $k$ be a smooth vector field.
Let $\{H_i\}$ be a family of Poisson commuting  Hamiltonian functions on a symplectic manifold. 
 $$dz_t^\epsilon = \f 1 {\sqrt \epsilon} \sum_{i=1}^n X_{H_i}(z_t^\epsilon)\circ dW_t^i+\f 1 \epsilon X_{H_0}(z_t^\epsilon)dt+ k(z_t^\epsilon)dt.$$ 
 Let $x_t^\epsilon=(H_1(z_t^\epsilon), \dots, H_1(z_t^\epsilon))$ and let $y_t^\epsilon$ denote the angle components. 
It was shown in \cite{averaging} that an averaging principle  holds  if $\{H_1, \dots, H_n\}$ forms a completely integrable systems and $x_t^\epsilon$ converging to the solution of and ODE. Furthermore if $k$ is also a Hamiltonian vector field, then $x^\epsilon$ converges on the scale $[0, \f 1 {\epsilon}]$ 
to a Markov process whose generator can be explicitly computed.

\item [(2)]  {\bf Stirring geodesics.} 
Let $H_0$ is a horizontal vector field on the orthonormal frame bundle $OM$ of a manifold $M$. The orthonormal frames  over $x\in M$E are the set of directions, i.e. linear maps from $\R^n$ to $T_xM$. Then the equation $\dot u_t=H_0(u_t)$ on $OM$ is the  equation for geodesics, its projection to the base manifold is a geodesic with unit speed, it solves a second order differential equation on $M$.
Let $\{A_k\}$ be a collection of skew symmetric $d\times d$ matrices where $d$ is the  dimension of the manifold. Let $A_k$ denote also the vertical fundamental fields obtained on $OM$ by rotating  an initial tangent vector in the direction of the exponential map of $A_k$. 
Consider the equation on $OM$:
$$du_t^\epsilon=\f 1 \epsilon H_0(u_t^\epsilon)dt+\f 1 \epsilon\sum_{k=1}^{n(n-1)\over 2} A_k(u_t^\epsilon)\circ dW_t^k.$$
These are large oscillatory perturbations on the the geodesic equation, considered on the scale  $[0, \f 1\epsilon]$.

Let $\pi$ denote the projection of a frame in $O_xM\define \L(\R^n, T_xM)$ to $x$.
We are only concerned with the projection of $u_t^\epsilon$ to the base manifold $M$ and we set $x_t^\epsilon=\pi(u_t^\epsilon)$.
Then our equation reduces to the following:
$$\dot x_t^\epsilon =\f 1 \epsilon x_t^\epsilon g_t^\epsilon$$ on the manifold $M$, where $g_t^\epsilon$ is a fast diffusion on $SO(d)$. If $\{A_k\}$ is an o.n.b. of
 ${\mathfrak so}(d)$, it was shown
in \cite{Li-geodesic} that  $\pi(u_t^\epsilon)$  converges as $\epsilon \to 0$. Theorem \ref{LLN} allow this theorem to be extended to
a set of $L=\{A_k\}$, which is not necessarily an o.n.b. of ${\mathfrak so}(n)$,  but the elements of $L$ and their Lie bracket has dimension $n$.

\item [(3)]   {\bf Perturbation to equi-variant diffusions.} Perturbed equi-variant SDEs on principal bundles were  studied  in \cite{Li-OM-1, conservation}. See also \cite{Elworthy-Li-Intertwining, ELL-equivariant, Elworthy-LeJan-Li-book-2} for the study of equivariant diffusions.
\item[(4)]  {\bf Inhomogensous scaling of Riemanian metric and collapsing of manifolds.} In \cite{Li-homogeneous}, a singularly perturbed model on a lie group $G$ with a compact subgroup $H$ was studied. 
Let $G$ be endowed with an ${\mathrm Ad}_H$-invariant left invariant Riemannian metric, this exists if $H$ is compact. 
Let $X_i$ be elements of the Lie algebra of $H$ which we denote by ${\mathfrak h}$. Let  $Y$ be an element of its orthogonal complement. We identify an element of the Lie algebra with the left 
invariant vector fields generated by it. Consider
 $$dg_t^\epsilon={1\over  \epsilon} \sum_{k=1}^p X_k(g_t^\epsilon)\circ dW_t^k+\f 1 \epsilon Y(g_t^\epsilon) \;dt.$$
Then, if $\{X_k\}$  and their brackets generates  ${\mathfrak h}$, it was shown in \cite{Li-homogeneous} that $\pi(g_t^\epsilon)$ converges to a diffusion process on the orbit space $G/H$.  
\end{enumerate} 
\end{example}
In these examples, a geometric reduction method is used, see also \cite{conservation,Ciccotti-Lelieve-Vanden-Eijnden}.

\subsection{ Diffusion homogenisation theory}
The models in [2] and [4] above  can be reduced to random ODEs with a random right hand sides. The effective dynamic theory falls into the diffusive homogensiation theory.
The first example also falls within this theory when the perturbation vector field $k$ is again a Hamiltonian vector field.
According to \cite{Hasminski-66}, the study of ODEs with a random right hand side was already considered in an article by Stratonovich in 1961.
See \cite{Komorowski-Landim-Olla, Mathieu-Piatnitski}.

A diffusive homogenisation theory, in its simplest form, is about a family of random ordinary differential equation of the form
$$\dot x_t^\epsilon=\f 1 {\sqrt \epsilon} V(x_t^\epsilon, y_t^ \epsilon)$$
where $V$ is a function with $\int V(x,y)\mu(dy)=0$ and $y_t$ is an ergodic Markov process or a stationary strong mixing stochastic process, and $y_\cdot^\epsilon$ is distributed as $y_{\f \cdot \epsilon}$, with stationary measure $\mu$. The  zero averaging assumption is about the `oscillatory  property' on the $y$ process and also restrictions on the functions $V$, without this  $x_t^\epsilon$ may blow up as $\epsilon\to 0$ even if $V$ is bounded and smooth.
 Assuming the square root scaling is the correct scaling for  a non-trivial effective limit, then we expect the effective dynamics to be a Markov process in which case we employ the martingale method.

Suppose that $y_t$ is a Markov process with generator $\L_0$, then $y^\epsilon_t$ is a Markov process with generator $\f 1\epsilon \L_0$ with invariant measure $\pi$. 
To proceed further let $F$ be a real valued function, then
$$F(x_t^\epsilon)=F(x_0^\epsilon) +\f 1 {\sqrt \epsilon}\sum_k \int_0^t (DF)_{x_s^\epsilon} \Big(V(x_s^\epsilon, y_s^ \epsilon)\Big) ds.$$
We take $x_0^\epsilon=x_0$ for simplicity.
Using a semi-martingale decomposition,
$$\int_0^t dF((V(x_s^\epsilon, y_s^ \epsilon))ds=M_t^\epsilon+A_t^\epsilon,$$
it remains to show that 
\begin{enumerate}
\item $(x_\cdot^\epsilon, \epsilon\in (0,1])$ is a relatively compact set of stochastic processes; 
\item any of its limit point  $\bar x_t$ solves a martingale problem with a generator $\bar \L$ which comes down to show that  conditioning on the past;
\item the drift part $A_t^\epsilon$ converges to the conditional process of $\int_0^t \bar \L  f(x_s) ds$.

\end{enumerate}If $y_t$ is a Markov process with generator $\L$, then the semi-martingale decomposition comes from solving the Poisson equation:
$$\L \beta(\cdot,y)=dF \circ V(\cdot, y).$$
This equation is solvable, for each $y$, precisely under the center condition on $V$. Of course we will need to assume some technical conditions, for example 
$\L_0$ is elliptic and the state space for $y$ is compact.
Then one applies the It\^o formula to obtain:
$$
{\begin{split}\beta(x_t^\epsilon, y_t^\epsilon)&=\beta(x_0, y_0)+ \f 1 {\sqrt \epsilon} \int_0^t D_x\beta \left( V(x_s^\epsilon, y_{\f s \epsilon})\right)ds+\f 1 \epsilon \int_0^t \L_y \beta  (x_s^\epsilon, y_{\f s \epsilon}) ds+ N_t^\epsilon\\
&=\beta(x_0, y_0)+ \f 1 {\sqrt \epsilon} \int_0^t D_x\beta \left( V(x_s^\epsilon, y_{\f s \epsilon})\right)ds+\f 1 \epsilon \int_0^t dF(V (x_s^\epsilon, y_{\f s \epsilon}) )ds+ N_t^\epsilon.
\end{split}} $$
 the subscript $y$ in $\L_y$ indicates applying the operator to the second variable $y$, and similarly $D_x\beta$ indicates differentiation in the first variable.
Also,  $N_t^\epsilon=\int_0^t D_y\beta (x_s^\epsilon, y^\epsilon_{s}) )d\{y_s^\epsilon\}$, where $\{y_s^\epsilon\}$ denotes the martingale part of $y^\epsilon$,  is a local martingale.
This means
\begin{align}\label{martingale}
F(x_t^\epsilon)=F(x_0^\epsilon) +\sqrt \epsilon\big( \beta(x_t^\epsilon, y_t^\epsilon)-\beta(x_0, y_0)\big) - 
\int_0^t D_x\beta \left( V(x_s^\epsilon, y_{\f s \epsilon})\right)ds -N_t^\epsilon.
\end{align}

This identity can now be used to show $\{x^\epsilon, \epsilon \in (0,1]\}$ is tight.  Let us examine the terms in the above equation.  We expect that that $\sqrt \epsilon\big( \beta(x_t^\epsilon, y_t^\epsilon)-\beta(x_0, y_0)\big)$ is negligible when $\epsilon$ is small and $\int_0^t D_x\beta \left( V(x_s^\epsilon, y_{\f s \epsilon})\right)ds $ converges if $x^\epsilon$ converges as a process.
Indeed if $V$ does not depend on $x$, the ergodic assumption will imply that  $\int_0^t D_x\beta \left( V(x_s^\epsilon, y_{\f s \epsilon})\right)ds$ converges to the spatial average $ \int V(y)\pi(dy)$. Let us define 
$$\bar \L=\int D_x \beta (V(x,y)) \pi(dy).$$
It remains to show that any limit $x_t$ satisfies that $F(x_t)-F(x_s) -\int \bar \L F(x_s) ds$ is a martingale, which can be obtained from (\ref{martingale}) by ergodic theorems.

If both the $x$ and $y$-variables taking values in a manifold and  $V(x,y)=\sum_{k=1}^NY_k(x)G_k(y)$ where $\bar G_k=\int G_k(y)\pi(dy)=0$,  a diffusive homogenisation theorem is proved in \cite{Li-limits},  from which we extract a simple version and state it below.
\begin{theorem}\cite{Li-limits}
Suppose that  $\L_0$ is a smooth and elliptic  operator on a smooth compact manifold $\CY$ with invariant measure $\pi$ and $G_k$ smooth functions with $\int_\CY  G_kd\pi=0$. 
Let $Y_i$ be smooth vector fields on $\R^N$, growing at mostly linearly, and have bounded derivatives of all order. 
Consider
$$\dot x_t^\epsilon=\f 1 {\sqrt \epsilon} \sum_{k=1}^NY_k(x_t^\epsilon) G_k(y_t^ \epsilon), \qquad \qquad x_0^\epsilon=x_0.$$
Then, the following statements hold.

\begin{enumerate}
\item [(1)]As $\epsilon\to 0$, the solution $x^\epsilon$ converges to a diffusion measure $\bar \mu$ on any bounded time intervals.   Furthermore for
every $r<\f 14$, their Wasserstein distance is bounded above by $\epsilon^r$, i.e.
 $$\sup_{t\le T}d (\hat P_{y^\epsilon_{t\over \epsilon}}, \bar \mu_t)\lesssim \epsilon^{\f 14-}.$$
 \item [(2)] Let $\beta_j$ be the solutions to $\L_0 \beta_j=G_j$,  $\overline{Y_i \beta_j} =\int_G (Y_i \beta_j) d\pi$, 
 and $P_t$ the Markov semigroup for 
  $$\bar \L=-\sum_{i,j=1}^m \overline{Y_i \beta_j}\, Y_i(Y_j f),$$
where $Y_iF=dF(Y_i)$. Then, for every $F\in BC^4$ and any $T$,
 $$\sup_{t\le T}\left|\E F(x_t^\epsilon) -P_tF(x_0)\right|
 \lesssim \epsilon  \sqrt{|\log \epsilon|} (1+|x_0|^2) \left(1+ |F|_{C^4}\right). $$
\end{enumerate}
\end{theorem}
%
\begin{remark}
In \cite{Li-limits} a more general theorem is proved without assuming the uniqueness of the invariant probability measures for the fast variables.
The compactness of the state space for $y_t$ is not important, what we really used is the exponential rate of convergence of the diffusion at time $t$ to the invariant measure in the total variation norm.
\end{remark}

\section{Appendix}

Here we give a proof for Lemma \ref{lemma-in-appendix}, which is completely analogue to that of \cite[Prop. 4.10] {Gehringer-Li-tagged},  where the theorem is stated for the fractional OU-process, but the proof is also valid for the processes specified in Lemma  \ref{lemma-in-appendix}. The proof is added here for reader's the convenience. We will however omit lengthy algebraic manipulations identical to that for the proof of  \cite[Prop. 4.10] {Gehringer-Li-tagged}, pointing out only where to find them.

{\bf Proof of Lemma  \ref{lemma-in-appendix}.} 
We write $U=G_i$ and $V=G_j$ and denote by $\mathcal{F}_k$ the filtration generated by the fractional Brownian motion defining our fOU process. 
For $k \in \N $, we define the $\F_k$-adapted processes:
\begin{align*}
I(k)&= \int_{k-1}^{k} U(y_s) ds, \quad J(k) = \int_{k-1}^{k} V(y_s) ds,\\
\hat{U}(k) &= \int_{k-1}^{\infty} \E \left[ U(y_s) | \mathcal{F}_k \right] ds, \quad  \hat{V}(k) = \int_{k-1}^{\infty} \E \left[ V(y_s) | \mathcal{F}_k \right] ds,\\
M_k &= \sum_{l=1}^k \hat{U}(l) - \E \left[ \hat{U}(l-1) | \mathcal{F}_{l-1} \right], \quad
N_k = \sum_{l=1}^k \hat{V}(l) - \E \left[ \hat{V}(l-1) | \mathcal{F}_{l-1} \right].
\end{align*}
In particular, $M_k$ and $N_k$ are $\mathcal{F}_k$ adapted $L^2$ martingales. There are the following useful identities. For $k \in \N $
\begin{align*}\label{mart-dif2}
&\hat U(k) = I(k)+ \E [\hat U(k+1) \, |\, \F_k],  \\
&M_{k+1}-M_k=I(k)+\hat U(k+1)-\hat U(k) \\
& \sum_{j=1}^k I(j) = \int_0^k U(y_r) dr =M_k-\hat U(k)+\hat U(1)-M_1,
\end{align*}
and similarly for $J,\hat{V}$ and $N$, where the function $U$ is replaced by $V$.

%

In a nutshell, the proof is as follows.
Combining Lemma  \ref{lemma2-area}   below, and using  $\E (U(y_s) V(y_r) ) =\E (U(y_{s-r})V(y_0))$, we see that
$$\begin{aligned}
&\epsilon \int_0^{\f t \epsilon} \int_0^s U(y_s) V(y_r) dr ds\\
&= \epsilon \sum_{k=1} ^{L}   (M_{k+1}-M_k)N_k 
+t   \int_0^{\infty}\E \left(U(y_v) V(y_0) \right)  du
+\err_1(\epsilon) + \err_2(\epsilon). \end{aligned}$$
The proof of  Lemma \ref{lemma-in-appendix} is then concluded  with Lemma \ref{lemma-6.15} and the identity
\begin{align*}
\left( \int_0^1 \! \!\int_0^s  + \int_{1}^{\infty}\! \!\int_0^1 \right) \E \left(U(y_s) V(y_r) \right) dr ds 
&=-\f 12\int_0^{\infty}  \! \! \int_u^{u-2}\E \left(U(y_v) V(y_0) \right) du dv\\
&=\int_0^{\infty}\E \left(U(y_v) V(y_0) \right)  dv.
\end{align*}

Henceforth in this section we set $L=L(\epsilon)=[\f t \epsilon]$.

\begin{lemma}\label{lemma2-area}
	There exists a function $\err_1(\epsilon)$, which converges to zero in probability as $\epsilon\to 0$,   such that
	\begin{equation}\label{area1} \begin{aligned}
	&\epsilon \int_{0}^{\f t \epsilon} \int_0^s  U(y_s) V(y_r) dr ds=\epsilon\sum _{k=1} ^{L} I(k) \sum_{l=1}^{k-1} J(l)+t \int _0^1 \int_0^s \E \left( U(y_s)  V(y_r)\right) \,dr  ds
	+\err_1(\epsilon)\end{aligned}
	\end{equation}
\end{lemma}

The proof for this follows exactly the same way as that of Lemma 4.13  in \cite{Gehringer-Li-tagged} where the proof is for the fractional Ornstein-Uhlenbeck process,  but the proof  used only the stationary property and the ergodicity of the process.

\begin{lemma}\label{lemma-6.15}
	The following converges in probability:
	$$\lim_{\epsilon\to 0} \left(  \epsilon\sum_{k=1} ^{L} I(k) \sum_{l=1}^{k-1} J(l)-  \epsilon  \sum_{k=1} ^{L}   (M_{k+1}-M_k)N_k\right)
	= t \;  \int_1^{\infty} \int_0^1 \E\left( U(y_s) V(y_r)\right) dr ds.$$
\end{lemma}

This lemma can be proved in the same way as that of   \cite[Lemma 4.14]{Gehringer-Li-tagged}, again that was for the fractional Ornstein-Uhlenbeck process. The idea is to use the above identities to 
write the summation on the left hand side as follows:
$$\epsilon \left( \sum_{k=1} ^{L} -I(k)\, \hat V(k) + \sum_{k=1} ^{L} I(k)(\hat V(1)-N_1)- \sum_{k=1} ^{L} (\hat U (k+1)-\hat U(k) ) N_k\right)=I_1^\epsilon+I_2^\epsilon +I_3^\epsilon.$$	
To the first  term we apply shift invariance and ergodic theorem to obtain 
\begin{equation}
I_1^\epsilon\to   -t\, \E [I(1) \hat V(1)]=(-t)  \E \left( \int_0^1 U(y_r) dr\int_0^{\infty} V(y_s) ds \right).
\end{equation}
Also, $$\begin{aligned}
	I_2^\epsilon=\E \left| \epsilon  \sum_{k=1}^{L } I(k)(\hat V(1)-N_1) \right|^2	\lesssim &\epsilon^2\, \E[\hat V(1)]^2  \int_0^{L }  \int_0^{L }  \E[ U(y_r)U(y_s)]\, ds \,dr \to 0,\end{aligned}$$
	Since $  \int_0^{L }  \int_0^{L }  \E[ U(y_r)U(y_s)]\, ds \,dr\sim \f t \epsilon$.
We then change the order of summation to obtain the following decomposition
	$$ \begin{aligned}
	-I_3^\epsilon
	&=-\epsilon\ \sum_ {j=1}^{L-1} (N_{j+1}-N_j)   \hat U(L+1)  +\epsilon \sum_{j=1}^{L-1} (N_{j+1}-N_j) \,\hat U(j+1)
	-\epsilon \left (\hat U(L+1) - \hat U(1)\right)  N_1\\
	&=J_1^\epsilon+J_2^\epsilon+J_3^\epsilon.
	\end{aligned}$$
	By Birkhoff's ergodic theorem,  $J_1^\epsilon\to 0$  a.s. and similarly,
\begin{equation}
 J_2^\epsilon \to  t\, \E\left( \hat U(2) (N_2-N_1) \right).
\end{equation}
Since $\hat U(j)$ is bounded in $L^2(\Omega)$ by the assumption (for the fractional OU-process, this assumption is proved to hold), 
	$ \left|J_3^\epsilon \right|_{L^2(\Omega)} \lesssim \epsilon \to 0$. 
	This concludes the limit of the left hand side to be
	 \begin{equation}\label{diff}
	\begin{split}
	 t \,\E \left[  \hat U(2) (N_2-N_1)  -  I(1) \hat V(1) \right],
	\end{split}
	\end{equation}
which we need to rewrite. Firstly,
 $$\begin{aligned} (N_2-N_1) =\hat V(2)-\hat V(1) + \int_0^1 V(y_s)\, ds,  \qquad I(1)=\hat U (1) -   \E [\hat U (2) | \F_1 ].\end{aligned}$$ 
Secondly,	\begin{align*}
	&\E \left(  \hat U(2) (N_2-N_1)  -  I(1) \hat V(1) \right)\\
	&=\int_1^{\infty} \int_0^1\E \left(  U(y_s) V(y_r)  \right) dr ds + \E \left(  \hat U (2) \left(\hat V(2) -\hat V(1)\right) - \left(\hat U (1) -   \E [\hat U (2) | \F_1 ]\right) \, \hat V(1)  \right).
	\end{align*}
	Since $\hat V(1)$ is $\F_1$ measurable and by the shift covariance of $\hat U(k) \hat V(k)$,
	\begin{align*}
	\E \left(  \hat U (2) \left(\hat V(2) -\hat V(1)\right) - \left(\hat U (1) -   \E [\hat U (2) | \F_1 ]\right) \, \hat V(1)  \right)=0.
	\end{align*}
  This concludes the proof of the lemma.

\newcommand{\etalchar}[1]{$^{#1}$}
\def\dbar{\leavevmode\hbox to 0pt{\hskip.2ex \accent"16\hss}d} \def\cprime{$'$}
\def\cprime{$'$} \def\cprime{$'$} \def\cprime{$'$} \def\cprime{$'$}
\def\cprime{$'$} \def\cprime{$'$} \def\cprime{$'$} \def\cprime{$'$}
\def\cprime{$'$}

\end{document}